\documentclass{amsart}
\usepackage{amssymb,amsmath,amsfonts,euscript,mathrsfs,amsthm}
\usepackage[dvips]{graphicx}
\usepackage[dvips]{color}

\input xy
\xyoption{all}
\CompileMatrices

\def\cU{{\mathscr U}}

\def\fh{{\mathfrak h}}

\def\fl{{\mathfrak l}}

\def\fs{{\mathfrak s}}

\def\End{\operatorname{End}}

\def\vac{|0\rangle}

\numberwithin{equation}{section}
\newtheorem{thm}{Theorem}[section]

\newtheorem{prop}[thm]{Proposition}
\newtheorem{lem}[thm]{Lemma}

\theoremstyle{definition}
\newtheorem{defn}[thm]{Definition}

\newtheorem{ex}[thm]{Example}

\theoremstyle{remark}

\newtheorem{rem}[thm]{Remark}

\DeclareMathOperator{\vq}{V_Q^{\mathrm{super}}}
\begin{document}
\title{SUSY Lattice vertex algebras}
\author{Reimundo Heluani}
\address{Department of Mathematics, UC Berkeley \\ Berkeley, CA 94720 \\
USA}
\email{heluani@math.berkeley.edu}
\thanks{R.~Heluani was supported by the Miller Institute for basic research in science}
\author{Victor G. Kac}
\address{Department of Mathematics, MIT \\ Cambridge, MA 02139 \\ USA}
\email{kac@math.mit.edu}
\thanks{V.~Kac was supported in part by NSF grant DMS-0501395}
\date{}
\begin{abstract}
We construct and study SUSY lattice vertex algebras. As a simple example, we
obtain the simple vertex algebra associated to the vertex algebra
$V_c(N3)$ of central charge $c=3/2$, as the SUSY lattice
vertex algebra associated to $\mathbb{Z}$ with bilinear form $(a,b) =
2ab$. 
\end{abstract}
\maketitle

\section{Introduction}\label{sec:introd}
In our recent paper \cite{heluani3} we developped a $\Lambda$-bracket
formalism for $N_K=N$ supersymmetric (SUSY) vertex algebras, which
greatly simplifies calculations with quantum superfields. Recall
\cite{kac:vertex} that a vertex algebra $V$ is called $N_K=N$
supersymmetric if there exist $N$ odd operators $S^i$, satisfying $[S^i,
S^j] = 2 \delta_{ij} T$, where $T$ is the translation operator, such that
\[ {[S^i}, Y(a,z)] = Y (S^i a, z) \]
for any state $a \in V$. The presence of the operators $S^i$ allows one
to introduce superfields, and in \cite{heluani3} we use the
$\Lambda$-bracket formalism to perform computations with them.

In the present note we study $N_K=1$ SUSY lattice vertex algebras
$V_Q^{\mathrm{super}}$, associated to an integral lattice $Q$. We give a
simple characterization of these SUSY vertex algebras, similar to that of
``ordinary'' lattice vertex algebas \cite{kac:vertex}, and prove that, as
ordinary vertex algebras, 
\begin{equation}
	\vq \simeq V_Q \otimes F \bigl( \Pi \fh \bigr),
	\label{eq:intro1}
\end{equation}
where $V_Q$ is the ``ordinary'' lattice vertex algebra, associated to
$Q$, and $F( \Pi \fh)$ is the vertex algebra of free fermions, based on
the space $\fh = \mathbb{C} \otimes_{\mathbb{Z}} Q$ with reversed parity.

We prove the existence of a canonical $N=1$ conformal structure for
any $\vq$. Moreover, for each skewsymmetric operator $A$ on $\fh$ such
that $A^2=1$, we construct on $\vq$ an $N=2$ conformal structure,
and for each triple of such operators $A^i$, $i=1,2,3$, such that $A^i
A^j = \sqrt{-1} \varepsilon_{ijk} A^k$ for $i\neq j$ ($\sigma$-matrices),
we construct on $\vq$ a ``little'' $N=4$ conformal structure.

We write down the characters and supercharacters for all irreducible
$\vq$-modules and show that they span an $SL_2(\mathbb{Z})$-invariant
finite-dimensional space of functions, holomorphic on the upper-half
plane.

One of the simplest examples of $\vq$, when $Q=\mathbb{Z}$ with the
bilinear form $(a,b) = 2 ab$, turns out to be isomorphic to the simple
vertex algebra associated to the $N=3$ superconfomral algebra
$V_{3/2}(N3)$ with central charge $c = 3/2$. This result elucidates the
classification on unitary representations of the $N=3$ superconformal algebra
at $c = 3/2$, obtained in \cite{ss87}. 

Note also that this result shows that the vertex algebra $V_{3/2}(N3)$ is
rational (and even semisimple). To the best of our knowledge, this is the
only known non-zero value of $c$ (i.e. $c = 3/2$), for which $V_c (N3)$
is a rational vertex algebra.

Throghout the paper, we use the notation and terminology of our paper
\cite{heluani3} without further notice.
\section{Existence and Uniqueness}
Let $Q$ be an integral lattice of rank $r$ with a non-degenerate symmetric bilinear form
$(\cdot, \cdot)$. Let $W = \Pi \mathbb{C} \otimes_{\mathbb{Z}} Q$ be the
associated complex vector space with reversed parity. Extend $(\cdot,
\cdot)$ to $W$ by bilinearity. We have the $N_K=1$ SUSY vertex algebra $V^1(W)$ as in \cite[Example
5.9]{heluani3}, that is, the vertex algebra generated by odd superfields
$\alpha \in W$ satisfying the Lambda-brackets:
\begin{equation}
  {[\alpha}_\Lambda \beta] = (\alpha, \beta) \chi.
  \label{eq:0a}
\end{equation}

Let $\mathbb{C}_\varepsilon[Q]$ be the group algebra of
$Q$, twisted by the $\mathbb{C}^\times$-valued $2$-cocycle $\varepsilon$,
that is the unital associative algebra generated by symbols $e^\alpha,
\alpha \in Q$,
with multiplication rules:
\begin{equation*}
  e^\alpha e^\beta = \varepsilon(\alpha, \beta) e^{\alpha + \beta}.
\end{equation*}
We will construct an $N_K=1$ SUSY vertex algebra structure on the
vector space $\vq:=V^1(W) \otimes_{\mathbb{C}}
\mathbb{C}_\varepsilon[Q]$ such that we have the Lambda brackets
(\ref{eq:0a}) and
\begin{equation}
 {[\alpha}_\Lambda e^\beta] = (\alpha,\beta) e^\beta.
 \label{eq:0b}
\end{equation}
 
Equation (\ref{eq:0b})
corresponds to the following.
For each $h \in W$, $j \in \mathbb{Z}$ and $J=0,1$; we have the
operators $h_{(j|J)}$ on $V^1(W)$. We extend these operators to
$\vq$ by:
\begin{equation*}
  \begin{aligned}
    h_{(j|J)} \left( s \otimes e^\alpha \right) &= h_{(j|J)} (s)
    \otimes e^\alpha, \qquad (j|J) \neq (0|0), \\
    h_{(0|0)} \left( s \otimes e^\alpha \right) &= (h, \alpha) s
    \otimes e^\alpha.
  \end{aligned}
\end{equation*}
Note that this preserves the commutation relations between the operators
$h_{(j|J)}$, therefore we still have
\begin{equation*}
  {[h}_\Lambda h'] = (h,h') \chi, \qquad h, h' \in W.
\end{equation*}

Let $e^\alpha$ denote the operator of multiplication by $1 \otimes
e^\alpha$ on $\vq$, then it is easy to check that
\begin{equation*}
  \begin{aligned}
    {[h_{(j|J)}}, e^\alpha] &= 0, \qquad (j|J) \neq (0|0),\\
    {[h_{(0|0)}}, e^\alpha] &= (h, \alpha) e^\alpha.
  \end{aligned}
\end{equation*}
From these commutation relations, and denoting by $\Gamma_\alpha$ the
super-field corresponding to $1 \otimes e^\alpha$ for $\alpha \in Q$, we obtain the Lambda bracket
as in (\ref{eq:0b}):
\begin{equation}
  {[h}_\Lambda \Gamma_\alpha] = (h,\alpha) \Gamma_\alpha.
  \label{eq:2}
\end{equation}
This in turn gives:
  \begin{equation*}
    [h_{(j|J)}, \Gamma_\alpha(Z)] = (-1)^J (h,\alpha) Z^{j|J} 
    \Gamma_{\alpha}(Z).
  \end{equation*} 
\begin{lem}
  This last formula implies
  \begin{multline*}
    \Gamma_\alpha(Z) = e^\alpha \left(1 + \sum_{j \in \mathbb{Z}} Z^{-1-j|1}
    \alpha_{(j|1)} \right) \exp \left( - \sum_{j < 0}
    \frac{Z^{-j|0}}{j} \alpha_{(j|0)} \right) \times \\ \times \exp \left(
    - \sum_{j > 0} \frac{Z^{-j|0}}{j} \alpha_{(j|0)} \right) A_\alpha(Z),
  \end{multline*}
  for some operators $A_\alpha(Z)$ on $\vq$ commuting with all $h_{(j|J)}$.
  \label{lem:1}
\end{lem}
\begin{proof}
  Let
  \begin{multline*}
    X_\alpha =  \exp \left(  \sum_{j < 0}
    \frac{Z^{-j|0}}{j} \alpha_{(j|0)} \right) \left(1 - \sum_{j \in \mathbb{Z}} Z^{-1-j|1}
    \alpha_{(j|1)} \right) (e^\alpha)^{-1} \Gamma_\alpha(Z)  \times \\
    \times \exp \left(
     \sum_{j > 0} \frac{Z^{-j|0}}{j} \alpha_{(j|0)} \right).
  \end{multline*}
  For $(j|J) = (0|0)$ we have
  \begin{equation*}
    \begin{aligned}
      {[h_{(0|0)}}, X_\alpha] &= \exp \left(  \sum_{j < 0}
    \frac{Z^{-j|0}}{j} \alpha_{(j|0)} \right) \left(1 - \sum_{j \in \mathbb{Z}} Z^{-1-j|1}
    \alpha_{(j|1)} \right) \times \\ & \quad \times  (- (h, \alpha)) (e^\alpha)^{-1}
    \Gamma_\alpha(Z) \exp \left(
     \sum_{j > 0} \frac{Z^{-j|0}}{j} \alpha_{(j|0)} \right) + \\ & \quad
     + \exp \left(  \sum_{j < 0}
    \frac{Z^{-j|0}}{j} \alpha_{(j|0)} \right) \left(1 - \sum_{j \in \mathbb{Z}} Z^{-1-j|1}
    \alpha_{(j|1)} \right) \times \\ & \quad \times (e^\alpha)^{-1}
    (h,\alpha) \Gamma_\alpha(Z)  \exp \left(
     \sum_{j > 0} \frac{Z^{-j|0}}{j} \alpha_{(j|0)} \right) = 0.
    \end{aligned}
  \end{equation*}
  The other cases are similar (somewhat long) computations. 
  \end{proof}

As a corollary, note that 
\begin{equation*}
  \forall h \in W, \quad h_{(0|0)} A_\alpha(Z) s \otimes e^\beta = (h,
  \beta) A_\alpha(Z) s \otimes e^\beta,
\end{equation*}
therefore we must have $A_\alpha(Z) s \otimes e^\beta = a_{\alpha,
\beta}(Z)
s \otimes e^\beta$ for some $a_{\alpha, \beta} \in \mathbb{C}(
(Z))$.
In particular we see that $A_\alpha(Z) \vac = \vac + \theta d_\alpha \vac +
O(Z)$, where $O(Z)$ denotes a power series which is a multiple of $z$, and
$d_\alpha$ is an odd constant if 
we use a
Grassmann algebra $L$ as our base ring for the vertex algebra $V_Q$
instead of $\mathbb{C}$.

Applying $\Gamma_\alpha(Z)$ to the vacuum vector we obtain then:
\begin{equation}
  1 \otimes e^\alpha + \theta \left(d_\alpha (1 \otimes e^\alpha) +
  \alpha_{(-1|1)} \otimes e^\alpha \right) + O(z). \label{eq:1.bc}
\end{equation}

Now recall from \cite[Thm 4.16 (3)]{heluani3} that on $\vq$ we must have
$Y(Sa,Z) = D_Z Y(a,Z)$, where $D_Z = (\partial_\theta + \theta
\partial_z)$, therefore we
obtain the identity:
\begin{equation*}
  D_Z \Gamma_\alpha(Z) = d_\alpha \Gamma_\alpha(Z) + :\alpha(Z)
  \Gamma_\alpha(Z):
\end{equation*}
Replacing $\Gamma_\alpha(Z)$ by its expression from Lemma \ref{lem:1} in
this last equation, we obtain a differential equation for $A_\alpha$.
Indeed, note that:
\begin{multline}
  \partial_\theta \Gamma_\alpha(Z) = e^\alpha \left( \sum_{j \in
  \mathbb{Z}} Z^{-1-j|0} \alpha_{(j|1)} \right) \exp \left( - \sum_{j < 0}
    \frac{Z^{-j|0}}{j} \alpha_{(j|0)} \right) \times \\ \times \exp \left(
    - \sum_{j > 0} \frac{Z^{-j|0}}{j} \alpha_{(j|0)} \right) A_\alpha(Z)
    + e^\alpha \left( 1 + \sum_{j \in
  \mathbb{Z}} Z^{-1-j|1} \alpha_{(j|1)} \right) \times \\ \times\exp \left( - \sum_{j < 0}
    \frac{Z^{-j|0}}{j} \alpha_{(j|0)} \right)  \exp \left(
    - \sum_{j > 0} \frac{Z^{-j|0}}{j} \alpha_{(j|0)} \right)
    \partial_\theta A_\alpha(Z),
    \label{eq:2a}
\end{multline}
and similarly
\begin{multline}
  \theta \partial_z \Gamma_\alpha(Z) = e^\alpha \left( \sum_{j < 0}
  Z^{-1-j|1} \alpha_{(j|0)} \right) \exp \left( - \sum_{j<0}
  \frac{Z^{-j|0}}{j} \alpha_{(j|0)} \right) \times \\ \times \exp \left(
  - \sum_{j > 0} \frac{Z^{-j|0}}{j} \alpha_{(j|0)} \right) A_\alpha(Z) +
  e^\alpha \exp \left( - \sum_{j <0} \frac{Z^{-j|0}}{j} \alpha_{(j|0)}
  \right) \times \\ \times  \left( \sum_{j > 0} Z^{-1-j|1} \alpha_{(j|0)} \right)
  \exp \left( -\sum_{j > 0} \frac{Z^{-j|0}}{j} \alpha_{(j|0)} \right)
  A_\alpha(Z) + \\ +  e^\alpha    \exp \left(  - \sum_{j < 0}
  \frac{Z^{-j|0}}{j} \alpha_{(j|0)} \right)  \exp \left( - \sum_{j > 0}
  \frac{Z^{-j|0}}{j} \alpha_{(j|0)} \right) \theta \partial_z
  A_\alpha(Z).
  \label{eq:2b}
\end{multline}
Combining these last two expressions and recalling the definition of
$\alpha_-(Z)$ in \cite[3.2.4]{heluani3} we obtain after canceling the
exponentials:
\begin{equation*}
  D_Z A_\alpha(Z) = d_\alpha A_\alpha(Z) + Z^{-1|1} \alpha_{(0|0)}
  A_\alpha(Z).
\end{equation*}
Multiplying both sides of this equation by $z^{-\alpha_{(0|0)}}$ and
defining $B_\alpha(Z) = A_\alpha(Z) z^{-\alpha_{(0|0)}}$, we obtain the
diferential equation:
\begin{equation*}
  D_Z B_\alpha(Z) = d_\alpha B_\alpha(Z).
\end{equation*}
All of its solutions are of the form $B_\alpha = C_\alpha (1 + \theta d_\alpha)$ for some even
constant $C_\alpha$. It follows from (\ref{eq:1.bc}) that $C_\alpha = 1$
for all $\alpha \in Q$, therefore we arrive at the expression:
\begin{multline}
  \Gamma_\alpha(Z) = e^\alpha 
    Z^{\alpha_{(0|0)}|0} \left(1 + Z^{0|1} d_\alpha +  \sum_{j \in \mathbb{Z}} Z^{-1-j|1}
    \alpha_{(j|1)} \right) \times \\ \times\exp \left( - \sum_{j < 0}
    \frac{Z^{-j|0}}{j} \alpha_{(j|0)} \right)  \exp \left(
    - \sum_{j > 0} \frac{Z^{-j|0}}{j} \alpha_{(j|0)} \right).
    \label{eq:finalform}
\end{multline}
\begin{rem}
  Note that the lattice has to be integral, otherwise this expression would not be a
  well defined superfield.
  \label{rem:3}
\end{rem}

\begin{rem}
	It is convenient to write
	$\Gamma_\alpha$ in a ``normally ordered'' form:
	\begin{multline}
		\Gamma_\alpha(Z) = e^\alpha 
    Z^{\alpha_{(0|0)}|0} \left( 1  + Z^{0|1} d_\alpha \right) \exp\left(\sum_{j < 0} Z^{-1-j|1}
    \alpha_{(j|1)} \right) \times \\ \times\exp \left( - \sum_{j < 0}
    \frac{Z^{-j|0}}{j} \alpha_{(j|0)} \right)  \exp \left(
    - \sum_{j > 0} \frac{Z^{-j|0}}{j} \alpha_{(j|0)} \right) \exp \left(
    \sum_{j \geq 0} Z^{-1-j|1} \alpha_{(j|1)}
    \right).
		\label{eq:4}
	\end{multline}
\label{rem:4}
\end{rem}
\begin{rem}
	Note that from (\ref{eq:2a}) and (\ref{eq:2b}) it follows that
	\begin{equation}
		(\partial_\theta - \theta \partial_z) \Gamma_\alpha(z,\theta) = d_\alpha
		\Gamma_\alpha (z,\theta) + :\alpha(z,-\theta)
		\Gamma_\alpha(z,\theta):.
		\label{eq:2c}
	\end{equation}
	\label{rem:5}
\end{rem}
We need to check the axioms of an $N_K=1$ SUSY vertex algebra.
 First let's check translation invariance. Note that we already
know 
\begin{equation*}
  S (1\otimes e^\alpha) = d_\alpha (1 \otimes e^\alpha) + \alpha_{(-1|1)}
  \otimes e^\alpha.
\end{equation*}
We must have then:
\begin{equation*}
  S (s \otimes e^\alpha) = S(s) \otimes e^\alpha + (-1)^{p(s)} (s \otimes
  1) \Bigl(d_\alpha (1 \otimes e^\alpha) + \alpha_{(-1|1)} \otimes e^\alpha \Bigr).
\end{equation*}

We need to check translation invariance for the fields $\Gamma_\alpha
(Z)$. For this we first note (recall $e^\alpha$ is the operator of
multiplication by $1 \otimes e^\alpha$ on $\vq$):
\begin{multline*}
  [S,e^\alpha] (s \otimes e^\beta) = \varepsilon(\alpha,\beta) S (s
  \otimes e^{\alpha + \beta}) -  \\ - e^\alpha \left( S(s) \otimes e^\beta +
  (-1)^{p(s)} (s \otimes 1) (d_\beta (1 \otimes e^\beta) + \beta_{(-1|1)}
  \otimes e^\beta)\right) = \\ = \varepsilon(\alpha, \beta) \left( S(s)
  \otimes e^{\alpha + \beta} + (-1)^{p(s)} (s \otimes 1) (d_{\alpha+ \beta} (1
  \otimes e^{\alpha+ \beta}) + (\alpha+\beta)_{(-1|1)} \otimes
  e^{\alpha+\beta})
  \right) - \\ - \varepsilon(\alpha, \beta) \left( S(s) \otimes
  e^{\alpha+\beta} +
  (-1)^{p(s)} (s \otimes 1) (d_\beta (1 \otimes e^{\alpha+\beta}) + \beta_{(-1|1)}
  \otimes e^{\alpha+\beta})\right) = \\ = \varepsilon(\alpha, \beta)
  (-1)^{p(s)} s  \left( d_{\alpha+\beta} - d_\beta + \alpha_{(-1|1)} \right)
  \otimes
  e^{\alpha+\beta}.
\end{multline*}
Letting $B_\alpha (e^\beta) := d_{\alpha+\beta} e^\beta$, we can
write this as (note that the sign $(-1)^{p(s)}$ is accounted for since $B_\alpha$
and $\alpha_{(-1|1)}$ are odd operators):
\begin{equation*}
  [S, e^\alpha] = e^\alpha (\alpha_{(-1|1)} + B_\alpha - B_0 ).
\end{equation*}
Note also that $[S, z^{\alpha_{(0|0)}}] = 0$. Indeed:
\begin{multline*}
  S z^{(\alpha, \beta)} s \otimes e^\beta - \\ -z^{\alpha_{(0|0)}} \left( S(s)
  \otimes e^\beta + (-1)^{p(s)} (s \otimes
  1) \Bigl(d_\beta (1 \otimes e^\beta) + \beta_{(-1|1)} \otimes
  e^\beta \Bigr)\right)=0.
\end{multline*}
Now we can compute $[S, \Gamma_\alpha(Z)]$. 
\begin{multline*}
  [S, \Gamma_\alpha(Z)] = e^\alpha \left( \alpha_{(-1|1)} + B_\alpha -
  B_0 \right) Z^{\alpha_{(0|0)}|0} \left( 1+Z^{0|1}d_\alpha \right) \left(1 +
  \sum_{j < 0} Z^{-1-j|1}
    \alpha_{(j|1)} \right) \times \\ \times\exp \left( - \sum_{j < 0}
    \frac{Z^{-j|0}}{j} \alpha_{(j|0)} \right)  \exp \left(
    - \sum_{j > 0} \frac{Z^{-j|0}}{j} \alpha_{(j|0)} \right) \left( 1 + \sum_{j
    \geq 0} Z^{-1-j|1} \alpha_{(j|1)}
    \right)  -\\ - e^\alpha
    Z^{\alpha_{(0|0)}|0} \left( 1 + Z^{0|1} d_\alpha \right) \left( \sum_{j < 0}
     Z^{-1-j|1} \alpha_{(j|0)} \right) \times \\ \times \exp \left( - \sum_{j < 0}
    \frac{Z^{-j|0}}{j} \alpha_{(j|0)} \right)  \exp \left(
    - \sum_{j > 0} \frac{Z^{-j|0}}{j} \alpha_{(j|0)} \right) \left( 1 + \sum_{j
    \geq 0} Z^{-1-j|1} \alpha_{(j|1)}
    \right) + \\ + e^\alpha Z^{\alpha_{(0|0)}|0} \left( 1 + Z^{0|1} d_\alpha
    \right) \left(1 +  \sum_{j < 0} Z^{-1-j|1}
    \alpha_{(j|1)} \right) \times \\ \times \left( \sum_{j < 0} Z^{-j|0}
    \alpha_{(j-1|1)} 
    \right) \exp \left( - \sum_{j < 0}
    \frac{Z^{-j|0}}{j} \alpha_{(j|0)} \right)  \exp \left(
    - \sum_{j > 0} \frac{Z^{-j|0}}{j} \alpha_{(j|0)} \right) \times \\ \times
    \left( 1 + \sum_{j \geq 0} Z^{-1-j|1}\alpha_{(j|1)} \right) +  e^\alpha
    Z^{\alpha_{(0|0)}|0} \left( 1+ Z^{0|1} d_\alpha \right) \left(1 +  \sum_{j <
    0}Z^{-1-j|1}
    \alpha_{(j|1)} \right) \times \\ \times \exp \left( - \sum_{j < 0}
    \frac{Z^{-j|0}}{j} \alpha_{(j|0)} \right) \left( \sum_{j > 0}
    Z^{-j|0} \alpha_{(j-1|1)} \right) \exp \left(
    - \sum_{j > 0} \frac{Z^{-j|0}}{j} \alpha_{(j|0)} \right) \times \\ \times
    \left( 1+ \sum_{j \geq 0} Z^{-1-j|1} \alpha_{(j|1)} \right) - e^\alpha
    Z^{\alpha_{(0|0)}|0} \left( 1+ Z^{0|1}d_\alpha \right) \left( 1 + \sum_{j <
    0} Z^{-1-j|1} \alpha_{(j|1)}
    \right) \times \\ \times \exp \left( - \sum_{j < 0} \frac{Z^{-j|0}}{j} \alpha_{(j|0}) \right)
    \exp \left( - \sum_{j > 0} \frac{Z^{-j|0}}{j} \alpha_{(j|0)} \right) \left(
    \sum_{j \geq 0} Z^{-1-j|1} \alpha_{(j|0)} \right),
\end{multline*}
where we used (4.6.1) of \cite{heluani3}.

We can commute the operators $B_\alpha$ to the left using 
\begin{equation*}
  e^\alpha B_\beta = B_{\beta - \alpha} e^\alpha,
\end{equation*}
to obtain:
\begin{equation*}
  [S, \Gamma_\alpha(z,\theta)] =  (B_0 - B_{-\alpha}) \Gamma_\alpha(z,\theta) +
  :\alpha(z,-\theta) \Gamma_\alpha(z,\theta):,
\end{equation*}
and comparing with (\ref{eq:2c}) we obtain that translation invariance for the
operators $\Gamma_\alpha(Z)$ holds if and only if 
\begin{equation*}
	B_0 - B_{-\alpha} = d_\alpha,
\end{equation*}
and this in turn holds iff $d_\alpha$ is additive in  $\alpha \in Q$. From now
on,
we will
assume for simplicity that our ring of scalars is $\mathbb{C}$, therefore
$d_\alpha = 0$ for all $\alpha \in Q$. 

We are left only to check locality between the operators $\Gamma_\alpha(Z)$. Note
first that
\begin{multline*}
	\left(1 + \sum_{j \geq 0} Z^{-1-j|1} \alpha_{(j|1)} \right) \left(1 + \sum_{j <
	0}
	W^{-1-j|1} \beta_{(j|1)}
	\right) = \\ = \left(1 + \sum_{j <
	0}
	W^{-1-j|1} \beta_{(j|1)}
	\right) \left(1 + \sum_{j \geq 0} Z^{-1-j|1} \alpha_{(j|1)} \right) \left(1 -
	i_{z,w} \frac{\theta \zeta (\alpha,\beta)}{ z-w} \right).
\end{multline*}
Using standard computations we find:
\begin{equation}
	\Gamma_\alpha(Z) \Gamma_\beta(W) = \varepsilon(\alpha,\beta) i_{z,w}
	(z-w)^{(\alpha,\beta)} \left( 1 - i_{z,w} \frac{\theta\zeta (\alpha,
	\beta)}{z-w} \right) \Gamma_{\alpha,\beta}(Z,W),
	\label{eq:8}
\end{equation}
where 
\begin{multline*}
	\Gamma_{\alpha,\beta}(Z,W) = e^{\alpha+\beta} Z^{\alpha_{(0|0)}|0}
	W^{\beta_{(0|0)}|0} \exp \left( \sum_{j < 0} Z^{-1-j|1} \alpha_{(j|1)} +
	W^{-1-j|1} \beta_{(j|1)}\right) \times \\ \times \exp \left( - \sum_{j < 0}
	\frac{Z^{-j|0}}{j}\alpha_{(j|0)} + \frac{W^{-j|0}}{j} \beta_{j|0} \right)
	\exp \left( - \sum_{j > 0} \frac{Z^{-j|0}}{j}\alpha_{(j|0)} +
	\frac{W^{-j|0}}{j} \beta_{(j|0)} \right) \times \\ \times \exp \left( \sum_{j \geq 0} Z^{-1-j|1}
	\alpha_{(j|1)} + W^{-1-j|1} \beta_{(j|1)} \right).
\end{multline*}
And finally note that we can rewrite (\ref{eq:8}) as:
\begin{equation}
	\Gamma_\alpha (Z) \Gamma_\beta (W) = \varepsilon(\alpha,\beta) i_{z,w}
	(z-w-\theta\zeta)^{(\alpha,\beta)} \Gamma_{\alpha,\beta}(Z,W),
	\label{eq:9}
\end{equation}
and
\begin{equation}
	\Gamma_\beta(W) \Gamma_\alpha(Z) = (-1)^{(\alpha,\beta)}
	\varepsilon(\beta,\alpha) i_{w,z} (z-w-\theta\zeta)^{(\alpha,\beta)}
	\Gamma_{\alpha,\beta}(Z,W).
	\label{eq:10}
\end{equation}
Therefore we are in the same position as in the non-SUSY case. Namely, we need to declare
the parity of $\Gamma_\alpha$ to be the parity of $(\alpha,\alpha)$, and locality
holds if and only if:
\begin{equation}
	\varepsilon(\alpha,\beta)=
	(-1)^{(\alpha,\beta)+(\alpha,\alpha)(\beta,\beta)}
	\varepsilon(\beta,\alpha).
	\label{eq:11}
\end{equation}
Using Lemma 4.2 in \cite{heluani3}, we find:
\begin{equation}
	[\Gamma_\alpha(Z),\Gamma_\beta(W)] = 
	\begin{cases}
		0, & \quad (\alpha, \beta) \geq 0\\
		\varepsilon(\alpha,\beta) \left( D_W^{(-1-(\alpha,\beta)|1)}
		\delta(Z,W) \right)
		\Gamma_{\alpha,\beta}(Z,W), & \quad (\alpha, \beta) < 0.
	\end{cases}
	\label{eq:12}
\end{equation}
\begin{ex}
	Let $Q = \mathbb{Z}$ with the generator $\alpha = 1$ and the usual
	bilinear form 
	$(\alpha, \alpha) = 1$. Let $V^{\mathrm{super}}_{\mathbb{Z}}$ be the corresponding
	$N_K=1$ SUSY lattice vertex algebra, where we put $\varepsilon= 1$. Then we have
	\begin{equation}
		\begin{aligned}
			{[\Gamma_\alpha(Z)}, \Gamma_{-\alpha}(W)] &= - (D_W \delta)
			\Gamma_{\alpha, -\alpha}(Z,W)  \\ &= (D_Z \delta)
			\Gamma_{\alpha, -\alpha}(Z,W) \\ &= D_Z ( \delta
			\Gamma_{\alpha, -\alpha} (Z,W)) + \delta D_Z
			\Gamma_{\alpha, -\alpha}(Z,W) \\ &= D_Z (\delta
			\Gamma_{\alpha, -\alpha} (W,W)) + \delta D_Z
			\Gamma_{\alpha, -\alpha} (Z, W)|_{Z=W} \\ &= - (D_W
			\delta) \Gamma_{\alpha, -\alpha}(W,W) + \delta D_Z
			\Gamma_{\alpha,-\alpha}(Z,W)|_{Z=W}.
		\end{aligned}
		\label{eq:13}
	\end{equation}
	and noting that $\Gamma_{\alpha,-\alpha}(W,W) = 1$, $D_Z \Gamma_{\alpha,
	-\alpha}(Z, W)|_{Z=W} = \alpha(W)$, we obtain:
	\begin{equation}
		[\Gamma_{\alpha} (Z), \Gamma_{-\alpha} (W)] = - D_W \delta(Z,W) +
		\delta(Z, W) \alpha(W).
		\label{eq:14}
	\end{equation}
	Therefore we have:
	\begin{equation}
		{[e^\alpha}_\Lambda e^{-\alpha}] = \alpha + \chi.
		\label{eq:15}
	\end{equation}
	Expanding the fields $\Gamma_{\pm \alpha}(Z) = \psi^\pm(z) +
	\theta \varphi^\pm(z)$, where $\psi^\pm$ are odd and $\varphi^\pm$ are
	even, and the field $\alpha(Z) = \psi^0(z) + \theta \alpha(z)$, where
	$\psi^0$ is odd and $\alpha$ is even, we obtain the following
	lambda brackets:
	\begin{equation}
		\begin{aligned}
			{[\psi^0}_\lambda \psi^0] &= 1, & [\alpha_\lambda \alpha]
			&= \lambda, & [\alpha_\lambda \psi^\pm] &= \pm \psi^\pm,
			\\
			{[\psi^0}_\lambda \psi^\pm] &= 0, & {[\psi^0}_\lambda
			\varphi^\pm] &= \pm \psi^\pm, & [\alpha_\lambda
			\varphi^\pm] &= \pm \varphi^\pm, \\
			{[\psi^+}_\lambda \psi^-] &= 1, & {[\psi^\pm}_\lambda
			\varphi^\mp] &= \pm \psi^0, & {[\varphi^+}_\lambda
			\varphi^-] &= \alpha + \lambda,
		\end{aligned}
		\label{eq:sl2super}
	\end{equation}
	and all other brackets are zero (except the ones given by skew-symmetry). We
	recognize these lambda brackets as the commutation relations of the
	generators of $V^2(\fs\fl_{2,\mathrm{super}})$.

	On the other hand, let us consider the algebra generated by the three
	fermions $\psi^0, \psi^\pm$ with
	commutation relations:
	\begin{equation}
		{[\psi^0}_\lambda \psi^0] = 1, \qquad {[\psi^+}_\lambda \psi^-] =
		1.
		\label{eq:17}
	\end{equation}
	This is the tensor product of a free fermion with the usual lattice
	vertex algebra $V_{\mathbb{Z}}$. Define the following fields:
	\begin{equation}
		\alpha = :\psi^+ \psi^-:, \quad \varphi^\pm = \pm :\psi^0
		\psi^\pm:
		\label{eq:16}
	\end{equation}
	It is straightforward to check that the commutation relations
	(\ref{eq:sl2super}) are satisfied by the fields $\alpha, \varphi^\pm,
	\psi^0$ and $\psi^\pm$. 

	Define the following fields in this algebra of three free fermions:
	\begin{equation*}
		\begin{aligned}
		G &= ::\psi^+ \psi^-: \psi^0:, \\
		2 L &= :(T\psi^0)\psi^0: + 
		:(:\psi^+ \psi^-:) (:\psi^+ \psi^-:) :.
		\label{eq:nsfermions}
	\end{aligned}
	\end{equation*}
	We claim that these fields 
	generate a Neveu Shwarz algebra of central charge $3/2$. Indeed this is the
	usual $N=1$ structure of the boson fermion system, for the boson $:\psi^+
	\psi^-:$ and the fermion $\psi^0$. We know that $\psi^0$ is primary of
	conformal weight $1/2$ with respect to this Virasoro, and finally, we
	find:
	\begin{equation*}
		\begin{aligned}
		{[\psi^+}_\lambda 2 L] &= - :\psi^+ : \psi^+ \psi^-: : - :
		:\psi^+ \psi^-: \psi^+: +  \int_0^\lambda \psi^+  \\
		&= (\lambda - T) \psi^+,
	\end{aligned}
		\label{eq:nsfermions2}
	\end{equation*}
	hence $\psi^+$ (and similarly $\psi^-$) are primary of conformal weight
	$1/2$. 
	
	As with any vertex algebra with a super conformal vector, we can
	construct a SUSY VA following \cite{heluani3}. To compute the
	corresponding superfields, we need to compute $[G_\lambda \psi^\pm]$ and
	$[G_\lambda \psi^0]$. We find easily $[\psi^0_\lambda G] = :\psi^+
	\psi^-:$ and ${[\psi^\pm}_\lambda G] = \mp :\psi^\pm \psi^0:$. Therefore
	the superfields are:
	\begin{equation*}
		\alpha(Z) := \psi^0(z) + \theta :\psi^+(z) \psi^-(z):, \qquad
		\Gamma_{\pm}(Z) = \psi^\pm(z) \pm \theta :\psi^0(z)
		\psi^\pm(z):,
	\end{equation*}
	and these super-fields satisfy the Lambda brackets of the superfields
	$\alpha(Z)$ and $\Gamma_{\pm \alpha}(Z)$ computed above.

%
	\label{ex:1}

\end{ex}

\section{Conformal structure}
Following \cite[Ex. 5.8]{heluani3} we define:
\begin{equation}
	G = \sum_i : (S \alpha_i) \alpha^i: 
	\label{eq:18}
\end{equation}
where $\{\alpha_i\}$ and $\{\alpha^i\}$ are dual bases of $W$. We know from
\cite[Ex. 5.8]{heluani3} that this is a Neveu-Schwarz vector with central charge
$c = \tfrac{3}{2} \mathrm{dim} W$ on $V^1(W)$, and that all superfields $h(Z), h \in W$ are primary of conformal weight
$1/2$. From (\ref{eq:2}) we obtain ${[\Gamma_\alpha}_\Lambda h ] = (-1)^{(\alpha,
\alpha)} (h, \alpha) \Gamma_\alpha$ (recall that $h$ is always odd in our
setting, while the parity of $\Gamma_\alpha$ is the parity of $(\alpha,\alpha)$). Therefore we can compute using the
non-commutative Wick formula:
\begin{equation*}
	\begin{aligned}
		{[\Gamma_\alpha}_\Lambda G] &= - \sum 
		(\alpha_i, \alpha) :\Bigl( (\chi + S)\Gamma_\alpha \Bigr)
		\alpha^i: + \sum (-1)^{(\alpha, \alpha)} 
		(\alpha^i, \alpha) :(S \alpha_i) \Gamma_\alpha: -  \\ & \quad - 
		\int_0^\Lambda (-1)^{(\alpha,\alpha)} (\eta - \chi) (\alpha_i,
		\alpha)(\alpha^i, \alpha) \Gamma_\alpha  \\
		&= -
		:(S\Gamma_\alpha) \alpha: - \chi :\Gamma_\alpha \alpha: +
		(-1)^{(\alpha,\alpha)} :(S\alpha) \Gamma_\alpha: - (-1)^{(\alpha,
		\alpha)} \lambda (\alpha, \alpha) \Gamma_\alpha ,
	\end{aligned}
	\label{eq:19}
\end{equation*}
and recalling that $:\alpha(Z) \Gamma_\alpha(Z): = D_Z \Gamma_\alpha(Z)$ we see
immediately that $:\alpha S\Gamma_\alpha: = 0$ (cf. 
Appendix \ref{ap:a}). Using
\begin{equation*}
	\begin{aligned}
	:\Gamma_\alpha \alpha: &= (-1)^{(\alpha,\alpha)} :\alpha \Gamma_\alpha:,
	\\
	:(S\Gamma_\alpha) \alpha: &= - (-1)^{(\alpha, \alpha)} :\alpha S
	\Gamma_\alpha: + \int_{-\nabla}^0 \chi (-1)^{(\alpha,\alpha)} (\alpha,
	\alpha) \Gamma_\alpha \\ &= - (-1)^{(\alpha,\alpha)} :\alpha
	S\Gamma_\alpha: + (-1)^{(\alpha,\alpha)} (\alpha, \alpha) T\Gamma_\alpha
	\\ &= (-1)^{(\alpha, \alpha)} (\alpha, \alpha) T \Gamma_\alpha,
\end{aligned}
\end{equation*}
we obtain
\begin{equation*}
	{[\Gamma_\alpha}_\Lambda G] = (-1)^{(\alpha,\alpha)} \left( :(S\alpha)
	\Gamma_\alpha:  - \chi S \Gamma_\alpha -
	(\alpha,\alpha) (\lambda + T) \Gamma_\alpha
	\right).
\end{equation*}
Therefore:
\begin{equation*}
	\begin{aligned}
	{[G}_\Lambda \Gamma_\alpha] &= :
	(S \alpha) \Gamma_\alpha: + :\alpha S\Gamma_\alpha: + (\chi + S)
	S\Gamma_\alpha + (\alpha,\alpha) \lambda \Gamma_\alpha \\
	&= \left(  2 T + (\alpha,\alpha) \lambda + \chi S \right) \Gamma_\alpha.
\end{aligned}
\end{equation*}
Hence $\Gamma_\alpha$ is a primary field of conformal weight
$(\alpha,\alpha)/2$ with respect to $G$. Note that this last equation
together with the fact that $G$ is a conformal vector on $V^1(W)$
imply that $G$ is a conformal vector on $\vq$. Indeed, the
commutation relation of $G$ with itself is not changed (as it involves
only fields from the subalgebra $V^1(W) \subset \vq$) and from
this last equation we can easily show that $G_{(0|1)} = S$ on $\vq$.

When $r$ is even, we can enlarge this conformal
structure to an $N=2$ structure as follows. Let $A$ be an
endormorphism of $W$ satisfying:
\begin{equation}
  (A \alpha, \alpha') = - (\alpha, A \alpha'), \qquad A^2 =
   \mathrm{Id},
  \label{eq:n2ext1}
\end{equation}
where $\mathrm{Id}$ is the identity operator in $W$. We can construct then an
even
super-field of $\vq$, primary of conformal weight $1$, given by
\begin{equation}
	J = \frac{1}{2} \sum_{i,j = 1}^r (\alpha_i, A \alpha^j) :\alpha_i
  \alpha^j:.
  \label{eq:n2ext2}
\end{equation}
\begin{prop}\hfil
	\begin{enumerate}
		\item The superfields $G$ and $J$ generate the $N=2$ conformal
			vertex algebra of
  central charge $\tfrac{3}{2} r$ viewed as an $N_K=1$ SUSY vertex
  algebra as in \cite[Ex.
  5.10]{heluani3}. 
  \item  The super vector space $\vq$ carries a conformal $N_K=2$ SUSY
	  vertex algebra structure of central charge $\tfrac{3}{2}r$,
	  namely, the vector $\tau = \sqrt{-1} J_{(-1|1)}\vac$ is an
	  $N_K=2$ conformal vector. 
	  \end{enumerate}
  \label{prop:n3ext}
\end{prop}
\begin{proof}
	\begin{enumerate}
		\item
  We already know that $G$ is a Neveu-Schwarz super-field of central
  charge $\tfrac{3}{2} r$ and that $J$ is a primary
  super-field of conformal weight $1$. In order to compute
  ${[J}_\Lambda J]$ we first compute (here we sum over repeated
  indexes):
  \begin{equation*}
      {[\alpha_k}_\Lambda J] = -  (\alpha_i, A \alpha^k)\chi \alpha_i,
  \end{equation*}
  therefore
  \begin{equation*}
    {[J}_\Lambda \alpha_k] =  (\alpha_i, A \alpha^k) (\chi + S)
    \alpha_i.
  \end{equation*}
  It follows then:
  \begin{equation}
\begin{aligned}
	{[J}_\Lambda J] &= \frac{1}{2} (\alpha_k, A \alpha^m) (\alpha_i, A\alpha^k)
	:\Bigl( (\chi + S) \alpha_i \Bigr) \alpha^m: - \\ & \quad
	\frac{1}{2} ( \alpha_k, A
  \alpha^m) (\alpha^j, A\alpha_m) :\alpha_k (\chi+S) \alpha^j: + \\ &
  \quad  \frac{1}{2} (\alpha_k, A \alpha^m) (\alpha_i, A \alpha^k) \int_0^\Lambda
  (\eta - \chi) (\alpha_i, \alpha^m) \eta d \Gamma, \\
  &=  (\alpha_k, A\alpha^m) (\alpha_i, A \alpha^k) : \Bigl( (\chi + S)
  \alpha_i \Bigr) \alpha^m: + \frac{1}{2} (\alpha_k, A \alpha^i) (\alpha_i, A
  \alpha^k) \lambda \chi, \\  &= :\bigl( (\chi+S) \alpha_i \bigr)
  \alpha^i: + \frac{r}{2} \lambda \chi
  = G + \frac{c}{3} \lambda \chi.
\end{aligned}
\label{eq:nosequepaso}
  \end{equation}
  And according to \cite[Ex. 5.10]{heluani3} these are the commutation
  relations of the $N=2$ super-vertex algebra (viewed as an $N_K=1$ SUSY
  vertex algebra).

  \item The superfields $G$ and $J$ allow us to construct a
	conformal $N_K=2$ SUSY vertex algebra structure (see
	\cite[Def. 5.6]{heluani3}) on the space
	$\vq$ as follows. We define the operators
	\begin{equation*}
		S^1 = G_{(0|1)}, \qquad S^2 = J_{(0|0)},
	\end{equation*}
	and it follows from the first part of the proposition that these operators
	satisfy:
	\begin{equation*}
		[S^i, S^j] = 2 \delta_{ij} T. 
	\end{equation*}
	With these operators now we can construct superfields for each $a
	\in \vq$ as 
	\begin{equation*}
		Y(a, z, \theta^1, \theta^2) = Y(a, z) + \theta^1 Y(S^1 a
		, z) + \theta^2 Y(S^2 a, z) + \theta^2 \theta^1 Y(a,z),
	\end{equation*}
	where $Y(a,z)$ is the usual field associated to $a$ when we view
	$\vq$ as an ``ordinary'' vertex algebra. Defining the vector
	$\tau \in \vq$ as
	\begin{equation*}
		\tau = \sqrt{-1} J_{(-1|1)}\vac,
	\end{equation*}
	we obtain easily now that this vector satisfies the properties of
	\cite[Def. 5.6]{heluani3}, namely, it is an $N_K=2$ conformal
	vector. 
	\end{enumerate}
\end{proof}

\begin{defn}
	Let $V$ be an $N_K=1$ conformal SUSY vertex algebra of central
	charge $c$. A
	\emph{little $N=4$ conformal structure of central charge $c$} on V consist of three
	superfields $\{J^i, i=1, 2, 3\}$, such that each pair $\{G,
	J^i\}$ defines an $N_K=2$ conformal structure of central charge
	$c$ on $V$ as in Prop.
	\ref{prop:n3ext} (2) and moreover, the fields $\{G, J^i\}$, $i=1,
	2, 3$, 
	satisfy the commutation relations of the $N=4$ vertex algebra of
	central charge $c$ as
	in \cite[Ex 5.10]{heluani3}.
\end{defn}
We obtain easily then:
\begin{prop} 
  Let $A^i$, $i = 1, 2, 3$, be endomorphisms of
  $W$ satisfying (\ref{eq:n2ext1}) and in addition
  \begin{equation*}
    A^i A^j = \sqrt{-1} \varepsilon_{ijk} A^k, \qquad i \neq j,
  \end{equation*}
  where $\varepsilon$ is the totally antisymmetric tensor. Define the
  superfields $J^i$, $i=1, 2, 3$ by (\ref{eq:n2ext2}) with $A$
  replaced by $A^i$. Then the superfields $G, J^i$, $i=1, \dots, 3$
  define a little $N=4$ conformal structure of central charge $c =
  \tfrac{3}{2} r$ on $\vq$. 
  \label{prop:n4ext}
\end{prop}

\begin{ex}
	It follows from the computation in Example \ref{ex:1} that when $(\alpha, \beta) = -1$ and
	$\alpha+ \beta \neq 0$ we have
	\begin{equation}
		{[\Gamma_\alpha}_\Lambda \Gamma_\beta] = \varepsilon(\alpha,
		\beta) \left(  :\alpha \Gamma_{\alpha + \beta}: + \chi
		\Gamma_{\alpha + \beta} \right).
		\label{eq:22}
	\end{equation}
	Suppose now that $(\alpha, \alpha) = 2$, then we see from (\ref{eq:12})
	that we have (we put $\varepsilon(\alpha, -\alpha) = 1$):
         \begin{equation}
		\begin{aligned}
			{[\Gamma_\alpha(Z)}, \Gamma_{-\alpha}(W)] &= -
			(D_W^{1|1} \delta(Z,W))
			\Gamma_{\alpha, -\alpha}(Z,W)  \\ &= - (D_Z^{1|1}
			\delta(Z,W))
			\Gamma_{\alpha, -\alpha}(Z,W) \\ &=  - D_Z^{1|1}
			( \delta(Z,W)
			\Gamma_{\alpha, -\alpha} (Z,W)) - (\partial_z
			\delta(Z,W))  D_Z
			\Gamma_{\alpha, -\alpha}(Z,W) + \\ & \quad +
			\left( D_Z \delta(Z,W) \right)
			\partial_z \Gamma_{\alpha, -\alpha} (Z,W)  -
			\delta(Z,W)
			D_Z^{1|1} \Gamma_{\alpha,-\alpha}(Z,W)  \\ &= 
			- D_Z^{1|1} (\delta(Z,W) \Gamma_{\alpha, -\alpha}(Z,W)) -
			\partial_z (\delta(Z,W) D_Z\Gamma_{\alpha,-\alpha}(Z,W)) +
			\\ & \quad +  \delta(Z,W)
			D_Z^{1|1}\Gamma_{\alpha, -\alpha}(Z,W) + D_Z
			(\delta(Z,W)
			\partial_z \Gamma_{\alpha, -\alpha}(Z,W)) \\ &= -
			D_Z^{1|1} (\delta(Z,W) \Gamma_{\alpha, -\alpha}(W,W)) -
			\partial_z (\delta(Z,W) D_Z \Gamma_{\alpha,
			-\alpha}(Z,W)|_{Z=W}) + \\ & \quad + \delta(Z,W)
			D_Z^{1|1}\Gamma_{\alpha,-\alpha}(Z,W)|_{Z=W} + D_Z
			(\delta(Z,W) \partial_z \Gamma_{\alpha,-\alpha} (Z,W)|_{Z=W})  \\ &= -
			(D_Z^{1|1} \delta(Z,W)) \Gamma_{\alpha, -\alpha}(W,W) -
			(\partial_z \delta(Z,W)) D_Z \Gamma_{\alpha,
			-\alpha}(Z,W)|_{Z=W} + \\ & \quad + \delta(Z,W)
			D_Z^{1|1}\Gamma_{\alpha,-\alpha}(Z,W)|_{Z=W} + (D_Z
			\delta(Z,W)) \partial_z \Gamma_{\alpha,-\alpha} (Z,W)|_{Z=W}  \\ &= -
			(D_W^{1|1} \delta(Z,W)) \Gamma_{\alpha, -\alpha}(W,W) +
			(\partial_w \delta(Z,W)) D_Z \Gamma_{\alpha,
			-\alpha}(Z,W)|_{Z=W} + \\ & \quad + \delta(Z,W)
			D_Z^{1|1}\Gamma_{\alpha,-\alpha}(Z,W)|_{Z=W} - (D_W
			\delta(Z,W)) \partial_z \Gamma_{\alpha,-\alpha} (Z,W)|_{Z=W}.
		\end{aligned}
	\end{equation}
	Note that 
	\begin{equation*}
		\begin{aligned}
	D_Z \Gamma_{\alpha, \beta}(Z,W)|_{Z=W} &= :\alpha(W)
	\Gamma_{\alpha + \beta} (W):, \\ \partial_z \Gamma_{\alpha,
	\beta}(Z,W)|_{Z=W} &= :(D_W \alpha(W)) \Gamma_{\alpha + \beta}(W): -
	:\alpha(W)  :\alpha(W) \Gamma_{\alpha + \beta} (W): : \\  D_Z^{1|1}
	\Gamma_{\alpha, \beta}(Z,W)|_{Z=W} &= :\partial_w \alpha(W)
	\Gamma_{\alpha+\beta} (W): + \\ & \quad + : \alpha(W) :(D_W \alpha(W)) \Gamma_{\alpha
	+\beta}(W): : - \\ & \quad : \alpha (W) :\alpha (W) :\alpha (W)
	\Gamma_{\alpha+\beta}(W): : :.
\end{aligned}
\end{equation*}
We obtain that for $(\alpha, \alpha) = 2$ we have:
\begin{equation}
	{[\Gamma_\alpha}_\Lambda \Gamma_{-\alpha}] = T\alpha + :\alpha S \alpha:
	+ \chi S\alpha + \lambda \alpha + \lambda \chi,
	\label{eq:lambdabracket2}
\end{equation}
where we have used the fact that $:\alpha \alpha: = 0$. 

Expanding the fields $\Gamma_{\pm \alpha}(Z) = e^\pm (z) + \theta
\psi^\pm(z)$ and $\alpha(Z) = \psi^0(z) + \theta h(z)$, we find the
commutation relations (cf. Appendix \ref{ap:a}):
\begin{equation}
  \begin{aligned}
    {[e^+}_\lambda e^-] &= h + \lambda, & {[e^\pm}_\lambda \psi^\mp] &= -
    :\psi^0 h: \mp \lambda \psi^0, \\
    {[\psi^+}_\lambda \psi^-] &= - T h - :hh: + :\psi^0 T \psi^0: - 2 \lambda h
    - \lambda^2, & [h_\lambda e^\pm] &=  \pm 2 e^\pm, \\
    [h_\lambda \psi^\pm] &= \pm 2 \psi^\pm, & {[\psi^0}_\lambda \psi^\pm]
    &= \pm 2 e^\pm, \\ {[\psi^0}_\lambda e^\pm] &= 0, & [h_\lambda h] &= 2
    \lambda, \\ {[\psi^0}_\lambda \psi^0] &= 2, & [h_\lambda \psi^0] &= 0.
  \end{aligned}
  \label{eq:longcommute}
\end{equation}
With respect to the natural conformal structure defined above, the field
$\psi^0$ has conformal weight $1/2$, the fields $h, e^\pm$ have conformal
weight $1$ and generate the current $\fs\fl_2$ algebra at level $1$,  commuting with the fermion $\psi^0$,  and the fields $\psi^\pm$ have conformal
weight $3/2$. 


   In order to remove the quadratic terms in the lambda-brackets we follow
   \cite[sec. 8.5]{kacwakimoto1} and define the fields:
   \begin{equation}
	   \begin{aligned}
		   J^0 &= h, & J^\pm &= e^\pm, & \Phi &= \frac{1}{\sqrt{8}}
		   \psi^0,
		   \\ \tilde{G}^0 &= - \frac{1}{\sqrt{8}} :\psi^0 h:,
		   & \tilde{L} &= \frac{1}{4} :hh: - \frac{1}{4} :\psi^0
		   T \psi^0:,
		   \\
		   \tilde{G}^- &= - \frac{1}{\sqrt{2}} \psi^-, & \tilde{G}^+ &=
		   \frac{1}{\sqrt{8}} \psi^+.
	   \end{aligned}
	   \label{eq:n3def}
   \end{equation}
   It follows that these fields satisfy the commutation relations of the $N=3$ super
   vertex algebra
    as defined in \cite[sec. 8.5]{kacwakimoto1}, the central charge is
   $\tilde{c} = 3/2$ and corresponds (in the example of \cite{kacwakimoto1}) to
   the Hamiltonian reduction of $osp(3|2)$ at level $k = - 3/4$. This is computed
   explicitly in Appendix \ref{ap:a}.
	\label{ex:2}
\end{ex}
\begin{ex}
	Consider the rank 2 lattice generated by vectors $\alpha^\pm$ such that
	$(\alpha^\pm, \alpha^\pm) = \pm 1$ and $(\alpha^\pm, \alpha^\mp) = 0$. We
	define $\varepsilon$ following \cite[5.5]{kac:vertex}
	\begin{equation*}
		\varepsilon(\alpha^\pm, \alpha^\pm) = \pm 1, \qquad
		\varepsilon(\alpha^\pm, \alpha^\mp) = \pm 1,
	\end{equation*}
	and extend by bimultiplicativity to $Q$.

	We
	have the operators $\Gamma_{\pm \alpha^\pm}$ satisfying the commutation
	relations 
	\begin{equation*}
		\begin{aligned}
			{[\Gamma_{\pm \alpha^+}}_\Lambda \Gamma_{\pm \alpha^+}]
			&= 0, & {[\Gamma_{\alpha^+}}_\Lambda \Gamma_{-\alpha^+}] &=
			\alpha^+ + \chi, \\
			{[\Gamma_{\pm \alpha^-}}_\Lambda \Gamma_{\mp \alpha^-}]
			&= 0, & {[\Gamma_{\alpha^-}}_\Lambda \Gamma_{\alpha^-}] &=
			- \left( \chi + \frac{1}{2} S \right) \Gamma_{2
			\alpha^-}.
		\end{aligned}
	\end{equation*}
	\label{ex:raro2}
\end{ex}

\begin{ex}
	In this example we let $W$ be the Cartan algebra of $\fs\fl(2|1)$ with
	its non-degenerate pairing. Namely, we will consider the rank 2 lattice
	generated by two elements $\alpha$ and $\beta$, and the bilinear form is such that
	$(\alpha, \alpha) = (\beta, \beta) = 0$ and $(\alpha, \beta) = (\beta,
	\alpha) = -1$. We let $\varepsilon(\alpha, \beta) = - \varepsilon(\beta,
	\alpha) = 1$. We have the corresponding vertex operators $\Gamma_{\pm
	\alpha}$ and $\Gamma_{\pm \beta}$ which are even, as well as the operator
	$\Gamma_{\alpha+\beta}$ which is also even. They satisfy the commutation
	relations:
	\begin{equation*}
		{[\Gamma_\alpha}_\Lambda \Gamma_{\beta}] =  \chi \Gamma_{\alpha
		+ \beta} + :\alpha \Gamma_{\alpha+ \beta}:. 
	\end{equation*}
	Note that $\Gamma_{\alpha + \beta}$ does not commute with
	itself. Indeed, it satisfies
	\begin{equation*}
		{[\Gamma_{\alpha+\beta}}_\Lambda \Gamma_{\alpha + \beta}] = -
		\frac{1}{2} \left(  TS  + \chi T + \lambda S
		+ 2 \lambda \chi) \right) \Gamma_{2(\alpha+ \beta)} 
		- :(\alpha + \beta) T\Gamma_{2 (\alpha + \beta)}:.
	\end{equation*}
	\label{ex:raro}
\end{ex}

\section{Representation Theory}
The representation theory of SUSY lattice vertex algebras parallels the non-SUSY
case, here we sketch a guideline:
\begin{itemize}
	\item The definition of a module over a SUSY VA is verbatim the
		definition in the usual case, namely a super vector space $M$ together
		with an association to any vector $a \in V$, a superfield
		$Y^M(a,Z)$ with values in $\End(M)$. Such that to the vacuum
		vector we associate the identity in $M$ and that $(j|J)$-th
		products are preserved. The notion of positive energy follows
		through, namely if the SUSY VA $V$ is conformal, we require
		the operator $L_0^M$ to act diagonally, with eigenvalues bounded
		below, and with finite dimensional eigenspaces.
	\item For a lattice $Q$ and its dual $Q^*$, we construct a $\mathbb{C}_\varepsilon [Q]$-module
		$\mathbb{C}_{\varepsilon^*} [Q^*]$ in the usual way, by defining
		\[ \varepsilon^*(\alpha, \mu + \beta) = \varepsilon(\alpha,
		\beta), \qquad \alpha, \beta \in Q, \mu \in Q^*. \]
	\item We can construct then a $\vq$-module $M$ as follows:  as a vector
		space we declare $M = V^1(W) \otimes_{\mathbb{C}}
		\mathbb{C}_{\varepsilon^*}[Q^*]$. The operators
		$h^M_{(n|j)}$ act, in the usual way, 
		on the first factor for $(n|j) \neq (0|0)$ and \[ h^M_{(0|0)} (s
		\otimes e^\lambda) = (h, \lambda) s \otimes e^\lambda, \qquad
		\lambda \in Q^*. \]
		The vertex operators $\Gamma^M_\alpha(Z)$ are defined as in
		(\ref{eq:finalform}), where $e^\alpha$ acts on
		$\mathbb{C}_{\varepsilon^*}[Q^*]$ and the operators
		$\alpha_{(n|j)}$ are replaced by the corresponding operators
		$\alpha_{(n|j)}^M$ (recall also that $d_\alpha = 0$). 
	\item As in the usual case, $M$ decomposes as a sum of irreducible
		modules \[V^{\mathrm{super}}_{\mu +
		Q} = V^1(W) \otimes_{\mathbb{C}} \mathbb{C}_{\varepsilon^*}
		[\mu + Q], \qquad \mu + Q \in Q^*/Q. \] 
		The proof of
		irreducibility reduces to analizing the action of $h^M_{(0|0)}$,
		namely, if $u = \sum s_i \otimes e^{\lambda_i}$ with pairwise
		distinct $\lambda_i$, belongs to a non-zero $V_{\mu + Q}$
		submodule $\cU$, we act diagonally by $h^M_{(0|0)}$ and since
		the Vandermonde matrix is invertible, we see that
		each $s_i \otimes e^{\lambda_i} \in \cU$. Now $V^1(W)$ is
		irreducible (recall it is the usual boson-fermion system)
		therefore $\vac \otimes e^{\lambda_i} \in \cU$. Applying
		$\Gamma_\beta^M(Z)$ we obtain that $(\vac \otimes e^{\lambda +
		\beta}) \in \cU$ for all $\beta \in Q$, and these generate $V_{\mu
		+ Q}$. 
	\item As a corollary, $\vq$ is simple. Given the examples above we see
		that $V^{\mathrm{super}}_{\mathbb{Z}}$ where $\mathbb{Z}$ is generated by $\alpha$
		with $(\alpha, \alpha) = m$ is the (irred. quotient of the) super affine algebra for
		$\hat{\fs\fl}_2$ at level $2$ when $m=1$, and it is the (irred.
		quotient of the) $N=3$ super vertex algebra at central charge
		$3/2$ when $m=2$. 
	\item Comparing the action of $h^M_{(0|0)}$ we see that all
		modules $V^{\mathrm{super}}_{\mu + Q}$ are non-isomorphic, and if $Q$ is a
		positive lattice we see that they are positive energy modules.
	\item The proof that $V_Q$ is rational and that irreducible
		representations are parametrized
		by $Q^*/Q$ now follows exactly like in the non-SUSY case.
\end{itemize}

\subsection{Character formulas}
Now let $Q$ be an even, integral lattice.
For a coset $\lambda \in Q^*/Q$, let $V_\lambda$ be the corresponding module over
$\vq$. If $Q$ is of rank $1$ then $V_\lambda$ has for basis monomials of the
form:
\begin{equation} \alpha_{(-j_m|0)}\dots \alpha_{(-j_1|0)} \alpha_{(-k_n-1/2|1)}
	\dots \alpha_{(-k_1 -1/2 |1)} \vac
	\otimes |\gamma>, \qquad \gamma \in \lambda + Q
\label{eq:generalbasis}
\end{equation}
where $1 \leq j_1 \leq \dots \leq j_{m}$ are integers, and $1/2 \leq k_1 < \dots < k_{n}$ are
half-integers. The general case is similar by fixing a basis $\{\alpha^i\}$ for
$W$. Let $\tau$ be a coordinate of the upper half plane, $q=e^{2 \pi i \tau}$,
let us assume that $z \in W$, and let $u$ be a complex
parameter. We define the full character of a module $M$ over $\vq$ to be
\begin{equation}
	\chi_M (\tau, z, u) = e^{2\pi i u} \mathrm{tr}_M e^{2 \pi i z_{(0|0)}}
q^{L_0^M - c/24},
\label{eq:chrdef}
\end{equation}
where $L_0^M$ is the energy operator (recall that  $\vq$ is a conformal
$N_K=1$ SUSY vertex algebra) and $c =
\tfrac{3 r}{2}$ is the central charge. Similarly, we define the
supercharacter of $M$ replacing the trace by the supertrace in (\ref{eq:chrdef}).
Let $p(j)$ be the number of partitions of the integer $j$ without
repeated odd parts. Similarly, let $sp(j)$ be the number of partitions of the
integer $j$, without repeated odd parts and an even number of odd parts,
minus the number of such partitions with odd number of odd parts. It is
easy to show that the generating functions for $p(j)$ and $sp(j)$ are
given by
\begin{equation*}
  \begin{aligned}
    \sum_{j = 1}^\infty q^j p(j) &= \prod_{j=1}^\infty
    \frac{(1+q^{2j-1})}{(1-q^{2j})}, \\
    \sum_{j = 1}^\infty q^j sp(j) &=  \prod_{j=1}^\infty
    \frac{(1 - q^{2j-1})}{(1 - q^{2j})}.
  \end{aligned}
  \label{eq:agreg}
\end{equation*}
Using this, we see that the characters $\chi_\lambda$ and the supercharacters
$\chi^s_\lambda$ corresponding to $V_\lambda$, $\lambda \in Q^*/Q$ are described
by:
 \begin{equation}
	 \begin{aligned}
		 \chi_\lambda (\tau, z, u) &= \left( \frac{\eta(\tau)}{\eta(2 \tau)
		 \eta(\tau/2)}\right)^{r} \Theta_{\lambda}^Q (\tau,
		 z, u), \\
		 \chi^s_\lambda (\tau, z, u) &= \left(
		 \frac{\eta(\tau/2)}{\eta(\tau)^2} \right)^{r}
		 \Theta_\lambda^Q (\tau, z, u),
	 \end{aligned}
	 \label{eq:charactersgeneral}
 \end{equation}
where $\Theta_\lambda^Q$ are the classical $\Theta$ functions defined as
\begin{equation*}
	\Theta_\lambda^Q(\tau, z,u) = e^{2 \pi i u} \sum_{\gamma \in \lambda + Q}
	q^{\frac{(\gamma,\gamma)}{2}} e^{2 \pi i (z, \gamma)},
\end{equation*}
and $\eta(\tau)$ is the Dedekind eta function:
\begin{equation*}
  \eta(\tau) = q^{1/24} \prod_{j = 1}^\infty (1 - q^j).
\end{equation*}
Recall that these functions satisfy the following \emph{modular transformation
properties} (we omit the superscript $Q$ where no confusion can arise):
\begin{equation*}
	\begin{aligned}
		\Theta_\lambda \left( -\frac{1}{\tau}, \frac{z}{\tau}, u -
		\frac{(z,z)}{2 \tau} \right) &= (-i \tau)^{r/2} |Q^*/Q|^{-1/2}
		\sum_{\lambda' \in Q^*/Q} e^{-2 \pi i (\lambda, \lambda')}
		\Theta_{\lambda'}(\tau, z, u), \\
		\Theta_\lambda (\tau+1, z, u) &= e^{i \pi (\lambda,\lambda)}
		\Theta_\lambda (\tau,z,u),
		\end{aligned}
\end{equation*}
and
\begin{equation*}
		\eta\left( -\frac{1}{\tau} \right) = (-i \tau)^{1/2} \eta(\tau),
		\qquad
		\eta(\tau+1) = e^{\pi i /12} \eta(\tau).
\end{equation*}
We note also:
\begin{equation*}
	\eta\left( \frac{\tau}{2} + \frac{1}{2} \right) = e^{\pi i/12}
	\frac{\eta(\tau)^3}{ \eta(2 \tau) \eta(\tau/2)}.
\end{equation*}

It follows that the characters (\ref{eq:charactersgeneral})
satisfy:
\begin{equation*}
	\begin{aligned}
		\chi_\lambda\left( -\frac{1}{\tau}, \frac{z}{\tau},u -
		\frac{(z,z)}{2\tau} \right) &= |Q^*/Q|^{-1/2} \sum_{\lambda' \in
		Q^*/Q} e^{-2\pi i (\lambda, \lambda')}
		\chi_{\lambda'}(\tau,z,u) \\
		\chi_\lambda (\tau + 1, z, u) &= e^{\pi i (\lambda,\lambda) - \pi
		i r /8} \chi^s_\lambda (\tau, z, u) \\
		\chi^s_\lambda (\tau + 1, z, u) &= e^{\pi i (\lambda,\lambda) - \pi
		i r/12} \chi_\lambda(\tau,z,u).
	\end{aligned}
\end{equation*}
%
To analyse the action of $S = (\tau \mapsto -1/\tau) \in SL(2, \mathbb{Z})$ on
the supercharacters, 
 we need to consider the
\emph{Ramond} sector for each of the modules $V_\lambda$. For this
consider the the automorphism $\sigma$ of $\vq$ given by
$v 
\mapsto (-1)^{p(v)} v$\footnote{Note that $\sigma$ is a vertex algebra automorphism
but it is not a SUSY vertex algebra automorphism.}. For each module $V_\lambda$,
let us call the corresponding $\sigma$-twisted module
$V^\mathrm{tw}_{\lambda}$. If $r = 1$, the bases for these modules are of the
form:
\begin{equation} \alpha_{-j_m|0}\dots \alpha_{-j_1|0} \alpha_{-k_n|1}
	\dots \alpha_{-k_1|1} \vac
\otimes |\gamma>, \qquad \gamma \in \lambda + Q,
\label{eq:basis3}
\end{equation}
where $1 \leq j_1 \leq \dots \leq j_{m}$ and $0 \leq k_1 < \dots < k_{n}$ are
integers (note that we supressed the parenthesis in the subscripts of the
creation operators since these are twisted fields). For lattices of general rank,
the basis is computed similarly. The monomial in
(\ref{eq:basis3}) has parity $(-1)^n$. To compute the energy, first we note that
for a free fermion system generated by 
odd fields $[\psi_\lambda \psi'] = (\psi,\psi')$, we consider the corresponding twisted
fields
\begin{equation*}
	 \psi^\mathrm{tw}(z) = \sum_{m \in \mathbb{Z}}
	\psi^\mathrm{tw}_m z^{-1/2-m},
\end{equation*}
with commutation relations
\begin{equation*}
	[\psi^\mathrm{tw}_m, {\psi'}^\mathrm{tw}_n] = (\psi,\psi') \delta_{m, -n}.
\end{equation*}
We have the twisted Virasoro $L^\mathrm{tw}(z) = \tfrac{1}{2}
\sum_i :h_i(z)h^i(z): - \tfrac{1}{2} :\psi_i(z)
\partial_z\psi^i(z):^\mathrm{tw}$, where $\{\psi_i\}$ and $\{\psi^i\}$ are dual
bases of $W$ and similarly $\{h_i\}$ and $\{h^i\}$ are dual bases of $\Pi W$. 
It is easy to show that this is equal to:
\begin{equation*}
	L^\mathrm{tw}(z) = \frac{1}{2} \sum_{i=1}^{r} :h_i(z) h^i(z): -
	\frac{1}{2} \sum_{i=1}^{r} :\psi_i^\mathrm{tw}(z) \partial_z
	{\psi^i}^\mathrm{tw}(z): -
	\frac{r}{16} z^{-2}.
\end{equation*}
It follows that the monomial (\ref{eq:basis3}) has
energy:
\[\frac{r}{16} + \frac{(\gamma,\gamma)}{2} + \sum_{i=1}^n k_i +
 \sum_{i=1}^m j_i.\]
The general even integral lattice case is calculated similarly.
It follows that the 
characters of these modules (clearly the supercharacters
vanish), are given by (we specialize the $\Theta$ function to $z=u=0$):
\begin{equation*}
	\chi_\lambda^{\mathrm{tw}}(\tau) = 2^{r} \frac{\eta(2\tau)}{\eta(\tau)^2}
	\Theta_\lambda(\tau),
\end{equation*}
the factor $2^{r}$ appears since all the creation operators
$\alpha_{0|1}$ preserve the energy.
We obtain therefore
\begin{equation*}
	\chi_\lambda^s (-\frac{1}{\tau}) =
	\frac{|Q^*/Q|^{-1/2}}{2^{r/2}} \sum_{\lambda' \in Q^*/Q}
	e^{-2 \pi i (\lambda, \lambda')}
	\chi_{\lambda'}^{\mathrm{tw}} (\tau).
\end{equation*}
Finally we find
\begin{equation*}
	\chi^{\mathrm{tw}}_\lambda (\tau + 1) =  e^{\pi i (\lambda, \lambda)}
	\chi^{\mathrm{tw}}_\lambda (\tau),
\end{equation*}
and therefore we arrive to the main theorem of this section:
\begin{thm}
  For each $\lambda \in Q^*/Q$, the linear span of the characters
  $\chi_\lambda(\tau)$, supercharacters $\chi_\lambda^s(\tau)$ and
  twisted characters $\chi^{\mathrm{tw}}_\lambda (\tau)$ is invariant
  under the action of $SL(2,\mathbb{Z})$. They satisfy:
  \begin{equation*}
    \begin{aligned}
   \chi_\lambda (\tau + 1) &= e^{\pi i (\lambda,\lambda) - \pi
		i r /8} \chi^s_\lambda (\tau), \\
		\chi^s_\lambda (\tau + 1) &= e^{\pi i (\lambda,\lambda) - \pi
		i r/12} \chi_\lambda(\tau),\\
		\chi^{\mathrm{tw}}_\lambda (\tau + 1) &=  e^{\pi i (\lambda, \lambda)}
	\chi^{\mathrm{tw}}_\lambda (\tau),
    \end{aligned}
  \end{equation*}
  and 
  \begin{equation*}
    \begin{aligned}
      \chi_\lambda \left( -\frac{1}{\tau}\right)
		 &= |Q^*/Q|^{-1/2} \sum_{\lambda' \in
		Q^*/Q} e^{-2\pi i (\lambda, \lambda')}
		\chi_{\lambda'}(\tau), \\
		\chi_\lambda^s (-\frac{1}{\tau}) &=
	\frac{|Q^*/Q|^{-1/2}}{(2)^{r/2}} \sum_{\lambda' \in Q^*/Q}
	e^{-2 \pi i (\lambda, \lambda')}
	\chi_{\lambda'}^{\mathrm{tw}} (\tau).
    \end{aligned}
  \end{equation*}
  \label{thm:modular1}
\end{thm}

\section{$\vq$ as a vertex algebra}
From the character formulas above, we guess that these SUSY lattice vertex
algebras are just tensor products of the usual (non-SUSY) lattice vertex algebras
and free fermions. To study this phenomena, 
recall that given a SUSY vertex algebra $V$, we can view $V$ as a  vertex
algebra by putting $\theta = 0$ in all superfields. In particular, the SUSY
vertex algebra associated to an integral lattice $Q$ is generated (as a usual
vertex algebra) by the fermions $\bar{h}(z) = h(z,0)$, the bosons $h(z) = (Sh)(z,0)$ for
$h \in W$, and the fields $\Gamma_\alpha(z) = \Gamma_\alpha(z,0)$ and
$\bar{\Gamma}_\alpha(z) = (S\Gamma_\alpha)(z,0)$, for $\alpha \in Q$. Note first
that putting $\theta=0$ in (\ref{eq:4}) we obtain the usual expression for the
vertex operator $\Gamma_\alpha(z)$ (see for example \cite[(5.5.9)]{kac:vertex}).
Putting $\theta=0$ in the identity \[S\Gamma_\alpha =
:\alpha \Gamma_\alpha:,\] we obtain that 
\begin{equation}
\bar{\Gamma}_\alpha(z) =
:\bar{h}(z)\Gamma_\alpha(z):.
\label{eq:aswq}
\end{equation}
In particular, if we let $V_Q$ be the
vertex algebra (as opposed to SUSY) associated to the integral lattice $Q$ and
$F(W)$ be the free fermions generated by odd fields $\bar{h}$, $h \in W$
satisfying
\[ {[\bar{h}}_\lambda \bar{h}'] = (h, h'), \]
we can construct a surjective morphism $V_Q \otimes F(W) \rightarrow \vq$.
Indeed this morphism maps $h \in V_Q \mapsto h \in \vq$, $\bar{h} \in F(W) \mapsto
\bar{h} \in \vq$ and $\Gamma_\alpha \in V_Q \mapsto \Gamma_\alpha \in \vq$. To
check that this morphism preserves lambda brackets we see that in $\vq$: $[h_\lambda
h'] = (h,h')\lambda$ and $[\bar{h}_\lambda \bar{h}'] = (h,h')$ which follow from
(\ref{eq:0a}). Finally, it follows from (\ref{eq:2}) that in $\vq$: $[\bar{h}_\lambda
\Gamma_\alpha] = 0$. 
 To check surjectivity, we only need to show that the generators of
$\vq$ lie in the image of this morphism. But we already know that $h, \bar{h}$
and $\Gamma_\alpha$ are in the image. It follows from (\ref{eq:aswq})
that $\bar{\Gamma}_\alpha$ is in the image of this morphism.
Since both algebras are simple, we obtain an isomorphism $V_Q \otimes F(W)
\simeq \vq$. 

As a corolary we see that $V_2(\fs\fl_{2,\mathrm{super}})$ is isomorphic the
tensor product of two charged fermions and one free fermion, and similarly
 we find that the simple $N=3$ super vertex algebra at central charge $c = 3/2$
 is isomorphic to the tensor product of the simple current algebra
 $V_1(\fs\fl_2)$ with a free fermion.

 \begin{rem}
	 In particular we see that the simple $V_c(N=3)$ super vertex algebra,
	 at central charge $c=3/2$ is a rational (even semisimple) vertex
	 algebra. To the authors knowledge, this is the only known
	 instance of $c$ for which $V_c(N=3)$ is rational and the list of
	 all of its modules and corresponding characters is known.
	 \label{rem:n3rational}
 \end{rem}

 We note also that taking the centralizer of the fermion $\Phi$ in the $N=3$
 vertex algebra, we find that the (irred. quotient of the) quantum Hamiltonian reduction of the super Lie
 algebra $\mathfrak{o}\fs\mathfrak{p}(3|2)$ as in \cite[Sec. 8.5]{kacwakimoto1} for $k = -3/4$ is
 just the simple vertex algebra $V_1(\fs\fl_2)$.

\appendix

\section{$N=3$ vertex algebra at central charge $c=3/2$} \label{ap:a}
In this section we compute explicitly the lambda-brackets obtained in Example
\ref{ex:2}. As a SUSY lattice vertex algebra, this is $V^{\mathrm{super}}_{\mathbb{Z}}$ where
$\mathbb{Z}$ is generated by $\alpha$ with $(\alpha,\alpha)=2$. We expand the
fields $\Gamma_{\pm \alpha}(Z)= e^\pm(z) + \theta\psi^\pm(z)$ and $\alpha(Z)=
\psi^0(z) + \theta h(z)$. It follows from (\ref{eq:2}) that
\begin{equation*}
	[\psi^0(z) + \theta h(z), e^\pm(w) + \zeta \psi^\pm(z)] = \pm 2 (\theta -
	\zeta) \delta(z,w) (e^\pm(w) + \zeta \psi^\pm(w)).
\end{equation*}
Collecting the terms with $\theta, \zeta$ and $\theta \zeta$ we see that this
implies:
\begin{subequations}
\begin{equation}
	\begin{aligned}
		{[}\psi^0(z), e^\pm(w)] &= 0, & [h(z), \psi^\pm(w)] &= \pm 2
		\delta(z,w) \psi^\pm(w), \\
		[\psi^0(z), \psi^\pm(w)] &= \pm 2 \delta(z,w) e^\pm(w), & [h(z),
		e^\pm(w)] &= \pm 2 \delta(z,w) e^\pm(w).
	\end{aligned}
	\label{eq:ap1}
\end{equation}
Similarly, it follows from $[\alpha_\Lambda \alpha] = 2\chi$ that
\begin{equation*}
	[\psi^0(z) + \theta h(z), \psi^0(w) + \zeta h(w)] = 2\delta(z,w) +
	2 \theta\zeta \partial_w \delta(z,w),
\end{equation*}
and collecting terms we obtain:
\begin{equation}
		[\psi^0(z), \psi^0(w)] = 2 \delta(z,w), \qquad [\psi^0(z), h(w)] =
		0, \qquad [h(z), h(w)] = 2 \partial_w \delta(z,w).
		\label{eq:ap2}
\end{equation}
Note also that from ${[\Gamma_{\pm\alpha}}_\Lambda \Gamma_{\pm \alpha}] = 0$ it
follows
\begin{equation}
	{[e^\pm(z)}, e^\pm(w)] = 0, \qquad [\psi^\pm(z), \psi^\pm(w)] = 0.
	\label{eq:ap2'}
\end{equation}
Finally, we obtain from (\ref{eq:lambdabracket2})
\begin{multline*}
	[e^+(z) + \theta \psi^+(z), e^-(w) + \zeta \psi^-(w)] = (\theta - \zeta)
	\delta(z,w) \bigl( \partial_w (\psi^0(w) + \zeta h(w)) + :\psi^0(w) h(w):
	+ \\ + \zeta :h(w)h(w): - \zeta :\psi^0(w) \partial_w \psi^0(w) \bigr) +
	\\ + (\delta(z,w) + \theta \zeta \partial_w \delta(z,w)) (h(w) + \zeta
	\partial_w \psi^0(w)) + \\ + (\theta - \zeta) \partial_w \delta(z,w)
	(\psi^0(w) + \zeta h(w)) + (\partial_w \delta(z,w) + \theta \zeta
	\partial_w^{2} \delta(z,w)),
\end{multline*}
and collecting terms as before:
\begin{equation}
	\begin{aligned}
		{[}e^+(z), e^-(w)] &= \delta(z,w) h(w) + \partial_w \delta(z,w), \\
		[\psi^+(z), e^-(w)] &= \delta(z,w) \Bigl( \partial_w \psi^0(w) +
		:\psi^0(w) h(w): \Bigr) + \bigl(\partial_w \delta(z,w)\bigr)
		\psi^0(w), \\
		[e^+(z), \psi^-(w)] &= - \delta(z,w) :\psi^0(w) h(w): - \left(
		\partial_w \delta(z,w) \right) \psi^0(w), \\
		[\psi^+(z), \psi^-(w)] &= - \delta(z,w) \left( \partial_w h(w) +
		:h(w)h(w): - :\psi^0(w)\partial_w\psi^0(w):  \right) \\ &\quad -
		(\partial_w \delta(z,w)) h(w) - (\partial_w \delta(z,w)) h(w) -
		\partial_w^2 \delta(z,w).
	\end{aligned}
\end{equation}
\label{eq:todas}
\end{subequations}
Equations (\ref{eq:todas}) easily imply equations (\ref{eq:longcommute}). 

Define now the fields $\tilde{G}^0, \tilde{G}^\pm$, $\tilde{L}$, $J^0$, $J^\pm$
and $\Phi$ as in (\ref{eq:n3def}). First we note that $\tilde{L}$ is the Virasoro
field of central charge $3/2$ given by the conformal structure (\ref{eq:18}), and we know already
that the field $\Phi$ is primary of conformal weight $1/2$, the fields $J^0,
J^\pm$ are primary of conformal weight $1$, and the fields $\tilde{G}^0$,
$\tilde{G}^\pm$ are primary of conformal weight $3/2$. 

It follows immediately from (\ref{eq:longcommute}):
\begin{subequations}
	\begin{gather}
		\begin{aligned}
			{[J^0}_\lambda J^\pm] &= \pm 2 J^\pm, & {[J^0}_\lambda
			J^0] &= 2 \lambda, & {[J^+}_\lambda J^-] &= J^0 +
			\lambda,\\
			{[J^+}_\lambda \tilde{G}^-] &= -2 \tilde{G}^0 + 2 \lambda
			\Phi, & {[J^-}_\lambda \tilde{G}^+] &= \tilde{G}^0 +
			\lambda \Phi, & {[\tilde{G}^\pm}_\lambda \tilde{G}^\pm]
			&= 0, \\
			{[\tilde{G}^+}_\lambda \tilde{G}^0] &=
			\frac{1}{4} (T + 2 \lambda) J^+, & {[\tilde{G}^+}_\lambda
			\Phi] &= \frac{1}{4} J^+, & {[\tilde{G}^-}_\lambda \Phi]
			&= \frac{1}{2} J^-, \\
	\end{aligned}\\
			 {[\tilde{G}^+}_\lambda \tilde{G}^-] = \tilde{L} +
			\frac{1}{4} (T + 2\lambda) J^0 + \frac{\lambda^2}{4}.
		\label{eq:total1}
	\end{gather}
	We have also:
	\begin{equation}
			{[J^0}_\lambda \tilde{G}^0] = [h_\lambda -
			\frac{1}{\sqrt{8}} :\psi^0 h: ] = -
			2 \lambda \frac{1}{\sqrt{8}} \psi^0 = - 2\lambda \Phi,
	\end{equation}
	and similarly
	\begin{equation}
		{[\Phi}_\lambda \tilde{G}^0] = - \frac{1}{8} {[\psi^0}_\lambda
		:\psi^0 h:] = - \frac{1}{4} h = -\frac{1}{4} J^0.
	\end{equation}
	On the other hand
	\begin{equation}
		\begin{aligned}
			{[\tilde{G}^+}_\lambda \tilde{G}^0] &=
			-\frac{1}{8} {[\psi^+}_\lambda :\psi^0 h:] \\
			&= -\frac{1}{8} \left( 2 :e^+ h: + 2 :\psi^0 \psi^+: + 2
			\int_0^\lambda {[e^+}_\mu h] d \mu\right) \\
			&= -\frac{1}{4} \left( :e^+ h: + :\psi^0 \psi^+: - 2
			\lambda e^+ \right)\\
			&= - \frac{1}{4} ( :h e^+: + :\psi^0 \psi^+: - 2(T +
			\lambda)e^+ ) \\
			&= \frac{1}{4} ( T + 2 \lambda) J^+ 
			- \frac{1}{4} :\psi^0 \psi^+: \\
			&=\frac{1}{4} (T + 2\lambda) J^+,
		\end{aligned}
	\end{equation}
	where we have used quasi-commutativity in the 4th line and the fact that
	$:h e^+: = T e^+$ and $:\psi^0 \psi^+: = 0$ in the last line. Both
	identities follow from $S\Gamma_\alpha = :\alpha \Gamma_\alpha:$ and
	$:\alpha \alpha: = 0$. Indeed
	\begin{equation*}
		:\alpha \alpha: = - :\alpha \alpha: + \int_{-\nabla}^0
		[\alpha_\Lambda \alpha] d \Lambda = - :\alpha \alpha:,
	\end{equation*}
	and therefore we obtain by quasi-associativity
	\begin{equation*}
		:\alpha S\Gamma_\alpha: = :\alpha :\alpha \Gamma_\alpha: : =
		::\alpha \alpha: \Gamma_\alpha: = 0.
	\end{equation*}
	It follows now
	\begin{multline*}
		\partial_ze^+(z) + \theta \partial_z\psi^+(z) = T\Gamma_\alpha(Z) = S^2\Gamma_\alpha(Z) = S :\alpha(Z)
		\Gamma_\alpha(Z): = \\ =
		:(S\alpha)(Z) \Gamma_\alpha(Z): = :h(z)e^+(z): + \theta \Bigl(
		:(\partial_z \psi^0(z)) e^+(z): + :h(z) \psi^+(z): \Bigr), 
	\end{multline*}
	from where $T e^+ = :h e^+:$. We also have 
	\begin{equation*}
		0 = :\alpha(Z) S\Gamma_\alpha(Z): = :\psi^0(z) \psi^+(z): +
		\theta \Bigl(:h(z)\psi^+(z): - :\psi^0(z) \partial_z e^+(z):
		\Bigr),
	\end{equation*}
	from where $:\psi^0\psi^+: = 0$.

	In a similar fashion:
	\begin{equation}
		\begin{aligned}
			{[\tilde{G}^-}_\lambda \tilde{G}^0] &=
			\frac{1}{4} {[\psi^-}_\lambda :\psi^0 h:]  \\ &=
			- \frac{1}{4} \left( 2 :e^- h: - 2 :\psi^0\psi^-: + 2
			\int_0^\lambda {[e^-}_\mu h] d\mu  \right) \\
			&= - \frac{1}{2} \left( :h e^-: + 2 Te^-
			+ 2 \lambda e^- \right) \\
			&= - \frac{1}{2} \left( T + 2\lambda \right) J^-.
		\end{aligned}
	\end{equation}
	We need to compute ${[\tilde{G}^0}_\lambda \tilde{G}^0]$, for this we
	compute first ${[\psi^0}_\lambda :\psi^0 h:]= 2 h$ and ${[h}_\lambda
	:\psi^0 h:] = 2 \lambda \psi^0$, from where
	\begin{equation}
		\begin{aligned}
		{[\tilde{G}^0}_\lambda \tilde{G}^0] &= \frac{1}{8} \left( 2
		:h h: - 2 :\psi^0 (\lambda + T)\psi^0: + 2 \int_0^\lambda
		[h_\mu h] d\mu 
		\right) \\
		&= \frac{1}{4} \left( :hh: - :\psi^0 T\psi^0: + \lambda^2 \right)
		\\&= \tilde{L} + \frac{\lambda^2}{4},
	\end{aligned}
	\end{equation}
	\label{total}
\end{subequations}
and according to \cite[pp. 41]{kacwakimoto1} equations (\ref{total}) are the
commutation relations for the generators of the $N=3$ super vertex algebra at
central charge $c=3/2$.

\end{document}